\documentclass{elsarticle}
\usepackage{mathrsfs}
\usepackage{amsfonts}
\usepackage{amsmath}
\usepackage{amssymb}
\usepackage{mathrsfs}

\usepackage{color}
\usepackage[dvipdfm]{hyperref}
\usepackage{graphics}
\newtheorem{theorem}{Theorem}[section]
\newtheorem{definition}{Definition}[section]
\newtheorem{lemma}[theorem]{Lemma}
\newtheorem{corollary}[theorem]{Corollary}
  \newtheorem{remark}[theorem]{Remark}

\usepackage{amssymb}
\usepackage{mathrsfs}
\usepackage{amsmath}
\usepackage{arydshln} % 添加虚线的宏包

\begin{document}

\title{The Laplacian polynomial of graphs
derived from regular graphs and applications \tnoteref{t1}}
\tnotetext[t1]{Partially supported by NNSFC (Nos.11471016,
11401004, 11171097, and 11371028), Anhui Provincial Natural
Science Foundation (No. 1408085QA03), Natural Science Foundation
of Anhui Province of China (No. KJ2013B105).}
\author[ahu,hbu]{Jia-Bao Liu}
\ead{liujiabaoad@163.com}
\author[ahu]{Xiang-Feng Pan\corref{cor1}}
\ead{xfpan@ustc.edu}
\author[ahu]{Fu-Tao Hu}
\ead{hufu@mail.ustc.edu.cn} \cortext[cor1]{Corresponding author.
Tel: +86-55163861313.}
\address[ahu]{School of Mathematical Sciences, Anhui University, Hefei, Anhui, 230601, P. R. China}
\address[hbu]{Department of Public Courses, Anhui Xinhua
University, Hefei, 230088, P. R. China}
\begin{abstract}
Let $R(G)$ be the graph obtained from $G$ by adding a new vertex
corresponding to each edge of $G$ and by joining each new vertex
to the end vertices of the corresponding edge. Let $RT(G)$ be the
graph obtained from $R(G)$ by adding a new edge corresponding to
every vertex of $G$, and by joining each new edge to every vertex
of $G$. In this paper, we determine the Laplacian polynomials of
$RT(G)$ of a regular graph $G$. Moreover, we derive formulae and
lower bounds of Kirchhoff index of the graphs. Finally we also
present the formulae for calculating the Kirchhoff index of some
special graphs as applications, which show the correction and
efficiency of the proposed results.
\end{abstract}


\date{}
\maketitle

\section{Introduction}
All graphs considered in this paper are simple and undirected. Let
$G = (V(G),E(G))$ be a graph with vertex set $V(E) = \{v_1, v_2,
\dots , v_n\}$ and edge set $E(G) = \{e_1, e_2, \dots , e_m\}$.
The adjacency matrix of $G$, denoted by $A(G)$, is the $n \times
n$ matrix whose $(i,j)$-entry is $1$ if $v_i$ and $v_j$ are
adjacent in $G$ and $0$ otherwise. Let $B(G)$ denote the adjacency
matrix and vertex-edge incidence matrix of $G$, which is the
$n\times m$ matrix whose $(i,j)$-entry is $1$ if $v_i$ is incident
to $e_j$ and $0$ otherwise. Denote $D(G)$ to be the diagonal
matrix with diagonal entries $d_G(v_1),d_G(v_2),\dots, d_G(v_n).$
The Laplacian matrix of $G$ defined as $L(G)=D(G)- A(G)$. The
Laplacian characteristic polynomial of $L(G)$, is defined as
$$\phi(L(G); x) = det(xI_n-L(G)),$$ or simply $\phi(L)$, where
$I_n$ is the identity matrix of size $n$, and its roots, denoted
by $0=\mu_1(G)\leq \mu_2(G)\leq \dots \leq \mu_n(G),$ are called
the Laplacian eigenvalues of $G$. The collection of eigenvalues of
$L(G)$ together with their multiplicities are called the
$L$-spectrum of $G$. Similar terminology will be used for $A(G).$
The adjacency characteristic polynomial of $G$, denoted by
$\varphi(A(G); x)$, is
 defined as $$\varphi(A(G); x) = det(xI_n-A(G)),$$the eigenvalues
of $A(G)$ are $\lambda_1(G)\geq \lambda_2(G)\geq \dots \geq
\lambda_n(G)$.  The collection of eigenvalues of $A(G)$ together
with their multiplicities are called the $A$-spectrum of $G$. For
other undefined notations and terminology from graph theory, the
readers may refer to~\cite{D. C1980} and the references therein.

 Klein and Randi\'c~\cite{Klein1993} introduced a new distance function named
resistance distance based on electrical network theory. The
resistance distance between vertices $i$ and $j$, denoted by
$r_{ij}$, is defined to be the effective electrical resistance
between them if each edge of $G$ is replaced by a unit
resistor~\cite{Klein1993}. The resistance distances attracted
extensive attention due to its wide applications in physics,
 chemistry, etc.~\cite{Feng2014,Fengl2014,FengG2014,Y2013,Zhang2013}.
   For more information on resistance distances
of graphs, the readers are referred to the recent papers
~\cite{Ya2014,YangJ2014,Liu2014}.

Until now, a large amount of graph operations such as the
Cartesian product, the Kronecker product, the corona and
neighborhood corona graphs have been introduced
in~\cite{Wang2012,Liu2013,Lu2013,McLeman2011,Gopalapillai2011}.
The following definition comes from~\cite{D. C1980} (see the
definition in p. 63 in~\cite{D. C1980}).

\begin{figure}[ht]
\center
  \includegraphics[width=\textwidth]{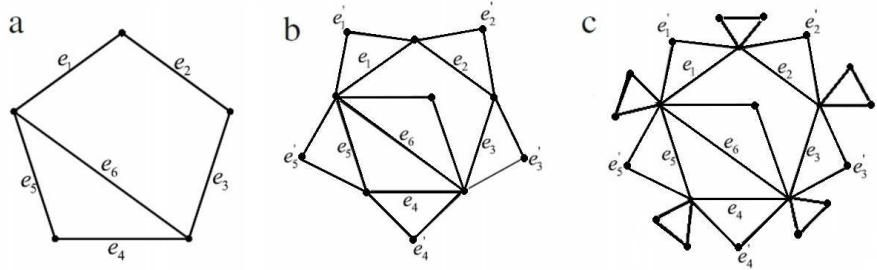}
  \vspace{-7em}
\caption{ (a) The graph $G$. ~~(b) The graph $R(G)$. ~~(c) The
graph $RT(G)$.}
\end{figure}

\begin{definition}\label{1-1}(see~\cite{D. C1980})
Let $R(G)=(V(R(G)),E(R(G)))$ be the graph obtained from $G$ by
adding a new vertex $e'$ corresponding to each edge $e = (a, b)$
of $G$ and by joining each new vertex $e'$ to the end vertices $a$
and $b$ of the corresponding edge $e=(a,b).$  (see Fig. 1(a) and
(b) for example).
\end{definition}

From the above definition, it is obvious that $R(G)$ is obtained
from $G$ by ``changing each edge $e=(a,b)$ of $G$ into a triangle
$ae'b$". Thus, $V(R(G))=V(G)\cup \{e'\mid e \in E(G)\}$ and
$E(R(G))=E(G)\cup\{(v_i, e'), (v_j, e')\mid e=(v_i, v_j)\in
E(G)\}$. A very elementary and natural question is what it would
 be like if we change each edge and each vertex of
$G$ into a triangle, which is stated as the following definiton.

\begin{definition}\label{1-2}
Let $RT(G)=(V(RT(G)),E(RT(G)))$ be the graph obtained from $RG$ by
adding a new edge $e''_i=(w^i_1, w^i_2)$ corresponding to each
vertex $v_i$ of $G$ and by joining the two vertices of each new
edge to each vertex $v_i$ of $G$, $RT(G)$ is obtained from $G$ by
``changing each edge and each vertex of $G$ into a triangle".
Thus, $V(RT(G))= \{e'\mid e \in E(G)\}\cup V(G)\cup \{w_1^i \mid
i=1, 2,\dots, n\} \cup \{w_2^i \mid i=1, 2,\dots, n\}$ and
$E(RT(G))=\{(v_i, e'), (v_j, e')\mid e=(v_i, v_j)\in E(G)\} \cup
E(G)\cup \{(v_i, w_1^i), (v_i, w_2^i)\mid v_i \in V(G),
e''_i=(w^i_1, w^i_2)\in E(RT(G)), i=1, 2,\dots, n\}$.(see Fig.
1(a), (b) and (c) for example).
\end{definition}

As the authors of~\cite{Zhang2009} pointed out, it is an
interesting problem to study the Kirchhoff index of graphs derived
from a single graph. In~\cite{XL2012}, the authors obtained
formulae and lower bounds of the Kirchhoff index of the line
graph, subdivision graph, total graph of a connected regular
graph, respectively.
 In~\cite{Wang2013}, Wang et al. determined the Laplacian polynomials
of $R(G)$ and $Q(G)$ of a regular graph $G$, they also derived
formulae and lower bounds of the Kirchhoff index of those graphs.
 Motivated by the
results, in this paper we further explore the Laplacian
polynomials of $RT(G)$ of a regular graph $G$. Moreover, we derive
the formulae and lower bounds of Kirchhoff index of the graphs. In
particular, special formulae are proposed for the Kirchhoff index
of $RT(G)$, where $G$ is a complete graph $K_n$, a cycle $C_n$ and
a regular complete bipartite graph $K_{n,n}$.

\section{Preliminaries}

At the beginning of this section, we review some concepts in
matrix theory. The Kronecker product $A\bigotimes B$ of two
matrices $A = (a_{ij})_{m\times n}$ and $B = (b_{ij})_{p\times q}$
is the $mp\times nq $ matrix obtained from $A$ by replacing each
element $a_{ij} $ by $a_{ij}B $. If $A, B, C$ and $D$ are matrices
of such size that one can form the matrix products $AC$ and $BD$,
then $(A\bigotimes B)(C\bigotimes D)=AC\bigotimes BD.$ It follows
that $A\bigotimes B$ is invertible if and only if $A$ and $B$ are
invertible, in which case the inverse is given by $(A\bigotimes
B)^{-1}=A^{-1}\bigotimes B^{-1}.$ Note also that $(A\bigotimes
B)^{T}=A^{T}\bigotimes B^{T}$. Moreover, if the matrices $A$ and
$B$ are of order $n\times n$ and $p\times p$, respectively, then
$det(A\bigotimes B) = (det~A)^p(det~B)^n.$ The readers are
referred to~\cite{Horn1991} for other properties of the Kronecker
product not mentioned here.

The symbols $\mathbf{0}_n$ and $\mathbf{1}_n$ (resp.,
$\mathbf{0}_{mn}$ and $\mathbf{1}_{mn}$) will stand for the
length-n column vectors (resp. $m\times n$ matrices) consisting
entirely of 0's and 1's.

\begin{lemma}\label{2-1}   (see~\cite{Zhang2005})  Let $M_1, M_2, M_3$ and $M_4$ be respectively $p\times p$, $p\times q$, $q\times p$ and $q\times q$ matrices with $M_1$ and $M_4$
invertible, then
\begin {eqnarray}
  det\left[%
\begin{array}{cc}
  M_1 & M_2 \\
  M_3 & M_4 \\
\end{array}%
\right]&&= det(M_4) \cdot det (M_1 - M_2M_4^{-1} M_3)\\
&&= det(M_1) \cdot det (M_4 - M_3M_1^{-1} M_2),
 \end {eqnarray} where $M_1 -
M_2M_4^{-1} M_3$ and $M_4-M_3M_1^{-1}M_2$ are called the Schur
complements of $M_4$ and $M_1$, respectively.
\end{lemma}

\section{\bf The Laplacian polynomials of RT(G)}

 For a regular graph $G$, the following
theorem gives the representation of the Laplacian polynomial of
$RT(G)$ by means of the characteristic polynomial and the
Laplacian polynomial of $G$, respectively.

 \begin{theorem}\label{3-1} Let $G$ be an $r$-regular graph with $n$ vertices
and $m$ edges, then $(i)~ \phi\Big(RT(G); \mu\Big)\\ =
(\mu-1)^n(\mu-2)^{m-n}(\mu-3)^n(3-\mu)^n\varphi\Big(G;\frac{(\mu-2)^2}{3-\mu}+\frac{r(2\mu-3)}{\mu-3}+\frac{2(\mu-2)}{(\mu-1)(\mu-3)}\Big).
$ $(ii)~ \phi\Big(RT(G); \mu\Big) =
(\mu-1)^n(\mu-2)^{m-n}(\mu-3)^{2n}\phi\Big(G;\frac{(\mu-2)^2}{\mu-3}-\frac{r\mu}{\mu-3}-\frac{2(\mu-2)}{(\mu-1)(\mu-3)}\Big).
$
\end{theorem}
{\bf Proof.} (i) Let $G$ be an arbitrary $r$-regular graph with
$n$ vertices and $m$ edges. Label the vertices of $RT(G)$ as
follows. Let $I(G) = \{e_1, e_2,\dots ,  e_{m}\}$, $V(G) = \{v_1,
v_2,\dots , v_{n}\}$  and $V(e'') = \{w_1, w_2\}$, and let $w^i_1,
w^i_2$ denote the vertices of the $i$-th copy of $e''$ for $i = 1,
2, \dots , n$, with the understanding that $w^i_j$ is the copy of
$w_j$ for each $j$. Denote $W_j =\{ w^1_j, w^2_j,\dots ,
w^{n}_j\}$, for $j = 1, 2$, then
\begin{equation}\label{}
    I(G) \bigcup V(G)  \bigcup
\big[W_1 \bigcup W_2 \big]
\end{equation}
 is a partition
of $V (RT(G))$. Obviously, the degrees of the vertices of $RT(G)$
are:
 $d_{RT(G)}(e_i) = 2 $, for $i = 1, 2,
\dots , m$,\\
 $~~~~~~~d_{RT(G)}(v_i) = 2d_{G}(v_i)+2$, for $i
= 1, 2, \dots , n$,\\  and $d_{RT(G)}(w_j^i) = 2$, for $i = 1, 2,
\dots , n, j = 1, 2$.

Let $B$ denotes vertex-edge incidence matrix of $G$.
 Since $G$ is an
$r$-regular graph, we have $D(G) = rI_{n} $. With respect to the
partition (3), then the Laplacian matrix of $RT(G)$ can be written
as

$$ L (RT(G))=  \left[%
\begin{array}{cc|lr}
  2I_{m} & -B^T & ~~~~~~~~0_{m\times 2n}  \\
  -B & L(G)+(r+2)I_{n} & ~~~~~　-I_{n}\otimes \mathbf{1}^T_{2} \\
   \hline
 ~~ 0_{2n\times m} & -I_{n}\otimes \mathbf{1}_{2}  & ~I_{n_1}\otimes \left[%
\begin{array}{cc}
  2 & -1 \\
  -1 & 2 \\
\end{array}%
\right]  \\
\end{array}%
\right],$$ where $\mathbf{1}_n$ denotes the all-one column vector
with size $n$.

By Lemma 2.1, we have
\begin {eqnarray}
\nonumber && \phi\Big(RT(G);\mu\Big)\\
\nonumber &&=det \left[%
\begin{array}{cc|lr}
  (\mu-2)I_{m} & B^T & ~~~~~~~~0_{m\times 2n}  \\\nonumber
  B & (\mu-r-2)I_{n}-L(G) & ~~~~~　I_{n}\otimes \mathbf{1}^T_{2} \\
   \hline
 ~~ 0_{2n\times m} & I_{n}\otimes \mathbf{1}_{2}  & ~I_{n}\otimes \left[%
\begin{array}{cc}
  \mu-2 & 1 \\
  1 & \mu-2 \\
\end{array}%
\right]  \\
\end{array}%
\right]\\
\nonumber &&=det \Bigg[I_{n_1}\otimes \left[%
\begin{array}{cc}
  \mu-2 & 1 \\
  1 & \mu-2 \\
\end{array}%
\right]  \Bigg]\cdot det ~S\\
 &&=(\mu-1)^n(\mu-3)^n \cdot det ~S,\\\nonumber
\end {eqnarray}
where
\begin {eqnarray*}
  && S \\
  &&=\left[%
\begin{array}{cc}
  (\mu-2)I_{m} & B^T \\
 B & (\mu-r-2)I_{n}-L(G) \\
\end{array}%
\right]
-\left[%
\begin{array}{c}
 \mathbf{0}_{0\times 2n}  \\
 I_{n}\otimes \mathbf{1}_{2}^T  \\
\end{array}%
\right] \\&&~~~ \cdot \left[%
\begin{array}{c}
  I_{n_1}\otimes \left[%
\begin{array}{cc}
  \mu-2 & 1 \\
  1 & \mu-2 \\
\end{array}%
\right]  \\
\end{array}%
\right]^{-1}
\left[%
\begin{array}{cc}
\mathbf{0}_{2n\times m}& I_n\otimes \mathbf{1}_{2}\\
\end{array}%
\right] \\
&&=\left[%
\begin{array}{cc}
  (\mu-2)I_{m} & B^T \\
 B & (\mu-r-2)I_{n}-L(G) \\
\end{array}%
\right]
-\left[%
\begin{array}{c}
 \mathbf{0}_{0\times 2n}  \\
 I_{n}\otimes \mathbf{1}_{2}^T  \\
\end{array}%
\right]\\
&&~~~\cdot   \left[%
\begin{array}{c}
  I_{n}\otimes \left[%
\begin{array}{cc}
  \mu-2 & 1 \\
  1 & \mu-2 \\
\end{array}%
\right]  \\
\end{array}^{-1}%
\right]
\left[%
\begin{array}{cc}
\mathbf{0}_{2n\times m}& I_n\otimes \mathbf{1}_{2}\\
\end{array}%
\right] \\
&&=\left[%
\begin{array}{cc}
  (\mu-2)I_{m} & B^T \\
 B & (\mu-r-2)I_{n}-L(G) \\
\end{array}%
\right]\\
&&~~-\left[%
\begin{array}{cc}
  \mathbf{0}_{m\times m} & \mathbf{0}_{m\times n} \\
\mathbf{0}_{n\times m} &  I_{n}\otimes \mathbf{1}_{2}^T\left[%
\begin{array}{cc}
  \mu-2 & 1 \\
  1 & \mu-2 \\
\end{array}%
\right]^{-1}\mathbf{1}_{2}\\
\end{array}%
\right]\\
&&=\left[%
\begin{array}{cc}
  (\mu-2)I_{m} & B^T \\
 B & (\mu-r-2)I_{n}-L(G) \\
\end{array}%
\right]
-\left[%
\begin{array}{cc}
  \mathbf{0}_{m\times m} & \mathbf{0}_{m\times n} \\
\mathbf{0}_{n\times m} &  I_{n}\otimes \frac
 {2}{\mu-1} \\
\end{array}%
\right]\\
&&=\left[%
\begin{array}{cc}
  (\mu-2)I_{m} & B^T \\
 B & (\mu-r-2-\frac
 {2}{\mu-1})I_{n}-L(G) \\
\end{array}%
\right].\\
\end {eqnarray*}

Let $l(G)$ be the line graph of $G$, it is
well-known~\cite{Nikiforov2007} that for a graph $G$,
$$BB^T = D(G)+A(G),  B^TB = 2I_m + A(l(G)).$$
Consequently,
\begin {eqnarray}
 \nonumber && det~ S\\
\nonumber &&=det \Big[ \big( \mu-2 )I_{m}\Big] \cdot det \Big[ ( \mu-r-2-\frac{2}{\mu-1} \big)I_{n}-L(G)-\frac{1}{\mu-2}BB^T )\Big]\\
 &&=(\mu-2)^m\cdot det \Big[ \big( \mu-r-2-\frac{2}{\mu-1}-\frac{(\mu-1)r}{\mu-2} \big)I_{n}+\frac{\mu-3}{\mu-2}A(G)  \Big]\\\nonumber
 &&=(\mu-2)^{m-n}(3-\mu)^n\\
 \nonumber &&~~\cdot det \Big[ \big( \frac{(\mu-2)^2}{3-\mu}+\frac{r(2\mu-3)}{\mu-3}+\frac{2(\mu-2)}{(\mu-1)(\mu-3)}\big)I_{n}-A(G)  \Big]\\
 \nonumber &&=(\mu-2)^{m-n}(3-\mu)^n\\
 &&~~\cdot \varphi\Big(G;\frac{(\mu-2)^2}{3-\mu}+\frac{r(2\mu-3)}{\mu-3}+\frac{2(\mu-2)}{(\mu-1)(\mu-3)}\Big).\\\nonumber
\end {eqnarray}
Actually, by virtue of (4) and (6) we have already established the
statement (i) in Theorem 3.1.

 (ii) Recall that $L(G)=rI_n-A(G)$.
It follows from (5) that
\begin {eqnarray}
\nonumber && det~ S\\
 \nonumber &&=(\mu-2)^m\cdot det \Big[ \big( \mu-r-2-\frac{2}{\mu-1}-\frac{(\mu-1)r}{\mu-2} \big)I_{n}+\frac{\mu-3}{\mu-2}A(G)  \Big]\\\nonumber
\nonumber &&=(\mu-2)^{m-n}(\mu-3)^n\\
\nonumber && ~~\cdot det \Big[ \big( \frac{(\mu-2)^2}{\mu-3}-\frac{r\mu}{\mu-3}-\frac{2(\mu-2)}{(\mu-1)(\mu-3)}\big)I_{n}-D(G)+A(G)  \Big]\\
\nonumber  &&=(\mu-2)^{m-n}(\mu-3)^n\\
\nonumber  &&~~\cdot det \Big[ \big( \frac{(\mu-2)^2}{\mu-3}-\frac{r\mu}{\mu-3}-\frac{2(\mu-2)}{(\mu-1)(\mu-3)}\big)I_{n}-L(G)  \Big]\\
 &&=(\mu-2)^{m-n}(\mu-3)^n\cdot \phi\Big(G;\frac{(\mu-2)^2}{\mu-3}-\frac{r\mu}{\mu-3}-\frac{2(\mu-2)}{(\mu-1)(\mu-3)}\Big).\\\nonumber
\end {eqnarray}
By combining (4) and (7), we get

$\phi\Big(RT(G); \mu\Big) =
(\mu-1)^n(\mu-2)^{m-n}(\mu-3)^{2n}\phi\Big(G;\frac{(\mu-2)^2}{\mu-3}-\frac{r\mu}{\mu-3}-\frac{2(\mu-2)}{(\mu-1)(\mu-3)}\Big).
$
 Thus the statement (ii) in Theorem 3.1 is proved.
 \hfill$\square$

\section{The Kirchhoff index of $RT(G)$}
In this section, we will explore the Kirchhoff index of the
$RT(G)$ of a regular graph $G$.

Zhu~\cite{ZKL}, Gutman and Mohar~\cite{GM} proved that the
relationship between Kirchhoff index of a graph and Laplacian
eigenvalues of the graph as follows.
\begin{lemma}\label{4-1}(\cite{ZKL,GM})
Let $G$ be a connected graph with $n\geq 2$ vertices, then
\begin{equation}
Kf(G)=n\sum^{n-1}_{i=1}\frac{1}{\mu_i}.
\end{equation}
\end{lemma}

 Denote by $\delta_i$ the degree of vertex $v_i\in V(G)$. Zhou
and Trinajsti\'c~\cite{BZ} proved that

\begin{lemma}\label{4-2}(\cite{BZ})
Let $G$ be a connected graph with $n\geq 2$ vertices, then
\begin{equation}
Kf(G)\geq -1+ (n-1) \sum_{v_i\in V(G)}\frac{1}{\delta_i},
\end{equation}
 with equality
attained if and only if $G = K_n$ or $G = K_{t, n-t}$ for $1\leq
t\leq\lfloor \frac{n}{2}\rfloor.$
\end{lemma}

The following lemma will be used later on.

\begin{lemma}\label{4-3} (\cite{XL2012})
Let $G$ be a connected graph with $n\geq2$ vertices and
$$\phi(G;\mu)=\mu^{n}+a_1\mu^{n-1}+a_2\mu^{n-2}+\dots
+a_{n-1}\mu, $$ then
$$\frac{Kf(G)}{n}=-\frac{a_{n-2}}{a_{n-1}}, ~~~(a_{n-2} = 1~ whenever ~n = 2),$$
where $a_{n-1},$ $ a_{n-2}$ are the coefficients of $\mu$ and $
\mu^2$ in the Laplacian characteristic polynomial, respectively.
\end{lemma}

Let $K_n $ be the complete graph with $n ~(n\geq 2)$ vertices. The
following theorem shows that $Kf(RT(G))$ can be completely
determined by the Kirchhoff index $Kf(G)$, the number of vertices
and the vertex degree of regular graph $G$.

\begin{theorem}\label{4-4}
 Let G be a connected $r$-regular graph with $n$ vertices, then
$$Kf(RT(G))=\frac{(r+6)^2}{6}Kf(G)+\frac{(r+5)n}{2}+\frac{(r+6)(5n-4)n}{6}+\frac{(r-2)(r+6)n^2}{8}.$$
 \end{theorem}

{\bf Proof.} Suppose first that $r = 1$, i.e. $G \cong K_2$.
 Since $Kf(RT(K_2))=\frac{74}{3}$.  It is easy to check that the result holds in
this case. Suppose now that $r\geq2$. Let
\begin{equation}\label{}
\phi(G;\mu)=\mu^{n}+a_1\mu^{n-1}+a_2\mu^{n-2}+\dots +a_{n-1}\mu.
\end{equation}
 It follows from Theorem 3.1 (ii) that
\begin {eqnarray}
\nonumber && \phi\Big(RT(G); \mu\Big)\\
 \nonumber &&=(\mu-1)^n(\mu-2)^{m-n}(\mu-3)^{2n}\phi\Big(G;\frac{(\mu-2)^2}{\mu-3}-\frac{r\mu}{\mu-3}-\frac{2(\mu-2)}{(\mu-1)(\mu-3)}\Big)\\
&&=(\mu-1)^n(\mu-2)^{m-n}(\mu-3)^{2n}\phi\Big(G;\frac{\mu\big[(\mu^2-(r+5)\mu+(r+6)\big]}{(\mu-1)(\mu-3)}\Big).\\\nonumber
\end {eqnarray}

Combining (10) with (11), one can obtain that
\begin {eqnarray*}
 &&\phi\Big(RT(G); \mu\Big)\\
&&=(\mu-1)^n(\mu-2)^{m-n}(\mu-3)^{2n}\Bigg\{\frac{\mu^n\big[\mu^2-(r+5)\mu+(r+6)\big]}{(\mu-1)^n(\mu-3)^n}^n  +\cdots \\
&&~~ +a_{n-2}\frac{\mu^2\big[\mu^2-(r+5)\mu+(r+6)\big]}{(\mu-1)^2(\mu-3)^2}^2\\
    &&~~  +a_{n-1}\frac{\mu\big[\mu^2-(r+5)\mu+(r+6)\big]}{(\mu-1)(\mu-3)}  \Bigg\}\\
&&=(\mu-2)^{m-n}(\mu-3)^{n}\Big\{{\mu^n\big[\mu^2-(r+5)\mu+(r+6)\big]}^n+\cdots \\
&&~~ +a_{n-2}\mu^2(\mu-1)^{n-2}(\mu-3)^{n-2}\big[\mu^2-(r+5)\mu+(r+6)\big]^2 \\
 &&~~ +a_{n-1}\mu(\mu-1)^{n-1}(\mu-3)^{n-1}\big[\mu^2-(r+5)\mu+(r+6)\big]  \Big\},\\
\end {eqnarray*}
where $\mu\neq 1, 3$. So the coefficient of $\mu^2$ in
$\phi(RT(G); \mu)$ is
\begin {eqnarray}
 \nonumber &&(-2)^{m-n}(-3)^{n}\Big[ a_{n-2}(r+6)^2(-1)^{n-2}(-3)^{n-2}\\
  \nonumber &&~~~~~~~~~~~~~~~~~~~~+ a_{n-1}(-r-5)(-1)^{n-1}(-3)^{n-1}\\
\nonumber &&~~~~~~~~~~~~~~~~~~~~+a_{n-1}(r+6)(n-1)(-1)^{n-2}(-3)^{n-1}\\
      \nonumber &&~~~~~~~~~~~~~~~~~~~~+a_{n-1}(r+6)(-1)^{n-1} (n-1) (-3)^{n-2}  \Big]\\
\nonumber &&~~~~~~~~~~~~~~~~~~~~ +(m-n) (-2)^{m-n-1}(-3)^{n} a_{n-1}(r+6)(-1)^{n-1}(-3)^{n-1} \\
&&~~~~~~~~~~~~~~~~~~~~+n(-3)^{n-1}(-2)^{m-n}a_{n-1}(r+6)(-1)^{n-1}(-3)^{n-1},\\\nonumber
\end {eqnarray}
and the coefficient of $\mu$ in $\phi(RT(G); \mu)$ is

\begin{equation}\label{}
    (-2)^{m-n} (-3)^na_{n-1}(r+6)(-1)^{n-1}(-3)^{n-1}.
\end{equation}

Notice that $RT(G)$ has $3n+m$ vertices. It follows from Lemma
4.3, (12) and (13) that
$$\frac{Kf(RT(G))}{3n+m}=-\frac{a_{n-2}}{a_{n-1}}\frac{r+6}{3}+\frac{r+5}{r+6}+\frac{5n-4}{3}+\frac{m-n}{2}.$$
Substituting the result of Lemma 4.3 and $m = \frac{nr}{2}$ into
the above equation.
$$\frac{Kf(RT(G))}{3n+\frac{nr}{2}}=\frac{r+6}{3}\frac{Kf(G)}{n}+\frac{r+5}{r+6}+\frac{5n-4}{3}+\frac{\frac{nr}{2}-n}{2}.$$
Simplifying the above result, one can obtain that
$$Kf(RT(G))=\frac{(r+6)^2}{6}Kf(G)+\frac{(r+5)n}{2}+\frac{(r+6)(5n-4)n}{6}+\frac{(r-2)(r+6)n^2}{8}.$$
Summing up, we complete the proof. \hfill$\square$

\begin{remark}\label{4-5}
Comparison to the Laplacian polynomials and its Kirchhoff indices
of $R(G)$ and $Q(G)$ in~\cite{Wang2013}, the graph $RT(G)$ has
more vertices and edges. It is clear that handling the problems of
Laplacian polynomial and Kirchhoff index are more difficult and
complex, but we deduce those with a simple approach.
 \end{remark}

In what follows, we propose a lower bound for the Kirchhoff index
for $RT(G)$ in terms of the number of vertices and the vertex
degree of a connected regular graph.

\begin{corollary}\label{4-6}
 Let G be a connected $r$-regular graph with $n$ vertices, then
$Kf(RT(G))\geq\frac{(r+6)^2(n^2-n-r)}{6r}+\frac{(r+5)n}{2}+\frac{(r+6)(5n-4)n}{6}+\frac{(r-2)(r+6)n^2}{8},$
and the equality holds if and only if $G \cong K_n$ or $G \cong
K_{\frac{n}{2}, \frac{n}{2}}$ and $n$ is even.
 \end{corollary}
{\bf Proof.} It follows from Lemma 4.2 and Theorem 4.4 that
\begin {eqnarray*}
&&Kf(RT(G))\\
&&\geq\frac{(r+6)^2}{6}\Big(\frac{(n-1)n}{r}
-1\Big)+\frac{(r+5)n}{2}+\frac{(r+6)(5n-4)n}{6}+\frac{(r-2)(r+6)n^2}{8}\\
&&=\frac{(r+6)^2(n^2-n-r)}{6r}+\frac{(r+5)n}{2}+\frac{(r+6)(5n-4)n}{6}+\frac{(r-2)(r+6)n^2}{8}.\\
\end {eqnarray*}
Clearly, the equality holds if and only if $G \cong K_n$ or $G
\cong K_{\frac{n}{2}, \frac{n}{2}}$ and $n$ is
even.~~\hfill$\square$

\section{Some applications}
In this section, we discuss some special graphs and give formulae
for their Kirchhoff index.
\subsection{Complete graph $K_n$ ($n\geq 2$)}
It is well known that $K_n$ is $(n-1)$-regular and $Kf(K_n) =n-1$.
It follows from Theorem 4.4 that
\begin {eqnarray*}
&& Kf(RT(K_n))\\
&&=\frac{(r+6)^2}{6}Kf(K_n)+\frac{(r+5)n}{2}+\frac{(r+6)(5n-4)n}{6}+\frac{(r-2)(r+6)n^2}{8}\\
&&=\frac{(r+6)^2(n-1)}{6}+\frac{(r+5)n}{2}+\frac{(r+6)(5n-4)n}{6}+\frac{(r-2)(r+6)n^2}{8}.\\
\end {eqnarray*}

Particularly, if $G\cong K_2$, one can obtain
$Kf(RT(K_2))=\frac{74}{3}$ by substituting $n=2, r=1$ into above
formula.

 In order to illustrate the correction and efficiency
of the above results, one can check $Kf(RT(K_2))$ for simplicity,
see Figure 2 (a).

It is easy to obtain

$$r_{12}=r_{13}=\frac{2}{3}, r_{14}=r_{15}=r_{16}=r_{17}=
\frac{4}{3}; r_{23}=r_{24}=r_{26}=\frac{2}{3}, r_{25}=r_{27}=
\frac{4}{3};$$
$$r_{34}=r_{36}=\frac{4}{3}, r_{35}=r_{37}=
\frac{2}{3};r_{45}=r_{47}=\frac{6}{3}, r_{46}=\frac{2}{3};
r_{56}=\frac{6}{3}, r_{57}= \frac{2}{3}; r_{67}= \frac{6}{3}.$$
Consequently, $Kf(RT(K_2))=\frac{74}{3},$ which coincides with the
above result.

\begin{figure}[ht]
\center
  \includegraphics[width=\textwidth]{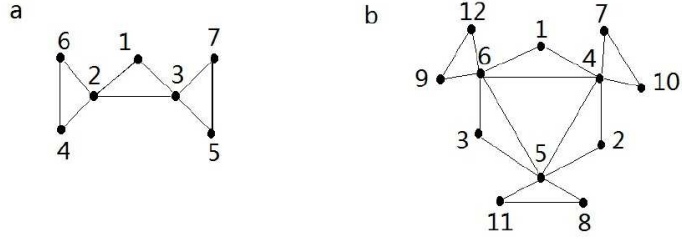}
  \vspace{-8em}
\caption{ (a) The graph $RT(K_2)$. ~~~(b) The graph $RT(C_3)$.}
\end{figure}

\subsection{Cycle $C_n$ ($n\geq 3$)}
 It was reported in~\cite{Lukovits} that $Kf(C_n)=
\frac{n^3-n}{12}.$ It follows from Theorem 4.4 that
\begin {eqnarray*}
&& Kf(RT(C_n))\\
&&=\frac{(r+6)^2}{6}Kf(C_n)+\frac{(r+5)n}{2}+\frac{(r+6)(5n-4)n}{6}+\frac{(r-2)(r+6)n^2}{8}\\
&&=\frac{(r+6)^2(n^3-n)}{72}+\frac{(r+5)n}{2}+\frac{(r+6)(5n-4)n}{6}+\frac{(r-2)(r+6)n^2}{8}.\\
\end {eqnarray*}

Similarly, for graph $RT(C_3)$, see Figure 2 (b). One can obtain

$Kf(RT(C_3))=\frac{455}{6},$ which also coincides with the above
formula.

\subsection{Complete bipartite graph $K_{n,n}$}
Note that  $K_{n,n}$ is $n$-regular with $2n$ vertices. Recall
from~\cite{XL2012} that
\begin{equation}\label{}
    Kf ( K_{n,n})=4n-3.
\end{equation}
It follows from (14) and Theorem 4.4 that
\begin {eqnarray*}
&&Kf(RT(K_{n,n}))\\
&&=\frac{(r+6)^2}{6}Kf(K_{n,n})+\frac{(r+5)\cdot 2n}{2}+\frac{(r+6)(10n-4)n}{3}\\
&&~~+\frac{(r-2)(r+6)n^2}{2}\\
&&=\frac{(r+6)^2(4n-3)}{6}+(r+5)n+\frac{(r+6)(10n-4)n}{3}+\frac{(r-2)(r+6)n^2}{2}.\\
\end {eqnarray*}

\section{Conclusions}
In this paper, based on the earlier definition $R(G)$, we
introduce a novel graph operation $RT(G)$, and explore its
Laplacian polynomial and Kirchhoff index.
%It is interesting that
%$Kf(RT(G))$ can be completely determined by the Kirchhoff index
%$Kf(G)$, the number of vertices and the vertex degree of regular
%graph $G$. 
By utilizing the spectral graph theory, we establish
the explicit formula for $Kf(RT(G))$ in terms of $Kf(G)$, the number of vertices and the vertex degree of regular
graph $G$, based on which we propose a
lower bound for the Kirchhoff index for $RT(G)$ with respect to the
number of vertices and the vertex degree.

\end{document}